\documentclass[12pt]{article}
\usepackage{amssymb}
\usepackage{amsthm}\usepackage{epsf,epsfig,amsfonts}
\usepackage{amsmath}
\usepackage{graphicx}

\usepackage{a4wide}

\newtheorem{corollary*}{Corollary}

\newcommand{\be}{\begin{equation}}
\newcommand{\ee}{\end{equation}}

\newcommand{\weg}[1]{}

\usepackage{epsf,epsfig,amsfonts,a4wide}

\usepackage{amsfonts}
\usepackage{amsmath,amssymb,amsthm,color}
\usepackage{url}

\newtheorem{Th}{Theorem}

\newtheorem{Lemma}{Lemma}
\newtheorem{Cor}{Corollary}

\theoremstyle{remark}

\newtheorem{Rem}{Remark}

\newcommand{\const}{\mbox{\rm const}}
\title{ Gallot-Tanno theorem for pseudo-Riemannian metrics and   a   proof that 
 decomposable cones over closed  complete pseudo-Riemannian manifolds  do not exist. }

\date{} \author{  Vladimir  S. Matveev\thanks{ Institute of Mathematics, FSU Jena, 07737 Jena Germany,  vladimir.matveev@uni-jena.de}}
\begin{document}
\maketitle

{\bf Introduction.} 
Let $g$ be a Riemannian or pseudo-Riemannian metric  
$$g= \sum_{i,j=1}^ng_{ij}(x_1,...,x_n)dx_i dx_j$$ on an $n-$dimensional manifold $M$. We consider the following equation  on the unknown function $\lambda$ on $M$.  \begin{equation} \label{tanno} 
  \nabla_k \nabla_j \nabla_i \lambda + 2 \nabla_k \lambda \cdot g_{ij} +
   \nabla_i \lambda \cdot g_{jk} + \nabla_j \lambda \cdot  g_{ik}=0.\end{equation}

This equation   is a famous one; it  naturally appeared in different parts  of differential geometry.  Couty \cite{couty} and   De Vries \cite{Vries}  studied it in the contex  of conformal transformations of  Riemannian metrics. They showed that, under certan additional assumptions,  conformal vector fields generate  nonconstant solutions  of the equation \eqref{tanno}.

    The equation also appears in investigation of geodesically equivalent metrics. Recall that two metrics on one manifold are \emph{geodesically equivalent},  if every geodesic of one metric is a reparametrized geodesic of the second metric.  Solodovnikov  \cite{sol} has shown that Riemannian metrics on $(n>3)-$dimensional manifolds admitting nontrivial   3-parameter family of geodesically equivalent metrics  allow nontrivial solutions of (a certain generlaization of) \eqref{tanno}.  
    Recently,  this result was generalised for pseudo-Riemannian metrics \cite[Corollary 4]{KM}. 
    Moreover,   as it was shown in \cite[Corollary 3]{einstein} (see also \cite{hall}), 
    an Einstein manifold of nonconstant scalar curvature 
    admitting nontrival geodesic  equivalence, after a proper scaling,  admits a nonconstant solution of 
    \eqref{tanno}.  Tanno  \cite{Tanno}  (see also \cite{hall})  related  the equation \eqref{tanno}  to {projective vector fields}, i.e., to vector fields whose local flows take unparametrized geodesics to geodesics.
    He has shown that  every nonconstant solution $\lambda $ of this equation   allows to construct    a nontrivial   projective vector field.

    Obata   used this equation trying to understand the relation  between the eigenvalues of the laplacian $\Delta_g$  and the geometry and topology of the manifold. 
   He  observed \cite{Obata} that the eigenfunctions corresponding to the second eigenvalue of the Laplacian of the  metrics of constant  curvature $+1$ on the sphere satisfy the equation \eqref{tanno}, and asked 
   the question whether the existence of a nonconstant solution of this equation 
   on a complete manifold implies that the manifold is covered by  the sphere with the standard metric. 
     The positive answer to  this question was indepedently and simultaneously   
   obtained by Gallot  \cite{G} and Tanno \cite{Tanno}.   
   
   This note    generalizes  the   result of Gallot  \cite{G} and Tanno \cite{Tanno} to pseudo-Riemannian metrics:

\begin{Th}  \label{main1}
Let $g$ be  
a light-line-complete 
connected  pseudo-Riemannian metric of indefinite signature (i.e., for no constant $c$ the metric $c\cdot g$ is Riemannian) on a closed   $n-$dimensional manifold $M^n$. Then, every solution of \eqref{tanno} is constant.  
\end{Th}

   \begin{Th} \label{main2}  Let $g$ be a negative-definite   metric (i.e., $-g$ is a 
   Riemannian  metric) on a   closed  connected manifold $M$. Then,  every solution of \eqref{tanno} is constant.  
   \end{Th} 
    
    Example  of Alexeevsky,  Cortes, Galaev and Leistner    \cite[Example 3.1]{ACGG} combined with Lemma \ref{cor2} below  shows that in the pseudo-Riemannian case the assumption that the metric is complete (but the manifold is  not closed)   is not sufficient to ensure that every  solution of \eqref{tanno} is  constant. 
   
   The equation \eqref{tanno} naturally appears  also in the investigation of the holonomy group  of      cones over pseudo-Riemannian manifolds.  Recall that the  
   \emph{cone over  $(M^n,g)$}  is the pseudo-Riemannian manifold  $(\hat M^{n+1}, \hat g)$, where $\hat M= \mathbb{R}_{>0}\times M$ and \begin{equation} \label{hatg} \hat g=(dx_0)^2 + x_0^2 \cdot\left(    \sum_{i,j=1}^ng_{ij}(x_1,...,x_n)dx_i dx_j\right),\end{equation}  where $x_0$ is the standard coordinate on  $\mathbb{R}_{>0}$ and $x_1,...,x_n$ are local coordinates on $M^n$.  
   Following \cite{ACGG,G},   we will show  that the decomposability of the cone (i.e., the existence of  a proper nondegenerate subspace  $U\subset T_p\hat M$ invariant with respect to the holonomy group) implies the existence of an nonconstant solution of \eqref{tanno} on $(M, g)$, see Lemma  \ref{cor2} below. Combining this with Theorems \ref{main1}, \ref{main2}, we obtain

\begin{Cor}  \label{Cor1} 
Let $g$ be  
a light-line-complete 
 pseudo-Riemannian metric of indefinite signature   on a closed  $n-$dimensional manifold $M^n$. 
 Then,    the corresponding cone  $(\hat M, \hat g)$ is not  decomposable.
\end{Cor}

\begin{Cor} \label{Cor2} 
Let $g$ be  
a complete  negative-definite 
 pseudo-Riemannian metric   on a closed   $n-$dimensional manifold $M^n$. 
 Then,    the corresponding cone  $(\hat M, \hat g)$ is not  decomposable.
\end{Cor} 

A partial case  of  Corollaries \ref{Cor1}, \ref{Cor2}  is   \cite[Theorem 6.1]{ACGG}. Our proof is  different from that of \cite{ACGG} and is shorter.  The case when the metric $g$ is Riemannian  was solved  in \cite[Proposition 3.1]{G}: Gallot used the Riemannian version  of Theorems \ref{main1}, \ref{main2} to show that if the cone $(\hat M, \hat g)$ over complete Riemannian $(M,g)$ is decomposable, then $g$ has constant curvature $+1$. \\[.1cm]

{\bf  Proof of Theorem \ref{main1}. }   Let $g$ be an  indefinite pseudo-Riemannian metric on $M^n $. Suppose the function  $\lambda$ satisfies \eqref{tanno}.  
 We take a light-line geodesic $\gamma(t)$ whose velocity vector will be denoted by $\dot\gamma= (\dot \gamma^i)$,  multiply  \eqref{tanno} by  $\dot \gamma^i \dot \gamma^j\dot \gamma^k$,  and sum over $i,j,k$. 
 Since the geodesic is light-line, at every point $\gamma(t)$ we have 
 $$\sum_{i,j=1}^n  
g_{ij} \dot \gamma^i \dot \gamma^j = \sum_{i,k=1}^n   g_{ik} \dot \gamma^i \dot \gamma^k=  \sum_{j,k=1}^n  g_{jk} \dot \gamma^j \dot \gamma^k \equiv 0 \textrm{ \ \ implying \ \ }  \sum_{i,j,k=1}^n \dot \gamma^i \dot \gamma^j\dot \gamma^k  \nabla_k\nabla_j \nabla_i \lambda  =0.$$   By definition of the geodesic, $\nabla_{\dot \gamma}{\dot \gamma}=0$ implying  
$$\sum_{i,j,k=1}^n \dot \gamma^i \dot \gamma^j\dot \gamma^k  \nabla_k\nabla_j \nabla_i \lambda  = \sum_{k=1}^n  \dot\gamma^k \nabla_k \left( \sum_{j=1}^n  \dot\gamma^j \nabla_j \left(\sum_{i=1}^n \dot\gamma^i \nabla_i \lambda  \right)  \right)  =      \tfrac{d^3}{dt^3} \lambda(\gamma(t))$$ implying $\tfrac{d^3}{dt^3} \lambda(\gamma(t))=0$ implying  that $\lambda=\const_2 t^2 + \const_1 t+ \const_0$. 

But by assumption the manifold $M$ is compact implying that the function $\lambda$ is bounded, 
 and the function $\const_2 t^2 + \const_1 t+ \const_0$ is bounded if and only if $\const_2=\const_1=0$. Then, $\lambda$ is constant along every light-line geodesic. Since every two points of a connected pseudo-Riemannian manifold of indefinite signature  can be connected by a sequence of light-line geodesics, the function $\lambda$ is   a constant. Theorem \ref{main1} is proved.  \\[.1cm]

{\bf Proof of Theorem 2.} We multiply \eqref{tanno} by $g^{ij}$ and sum over $i,j\in 1,...,n$.  We obtain: 
$
\nabla_k \left( \Delta_g \lambda \right) =-  2(n+1) \nabla_k \lambda, \textrm{ \ where $\Delta_g:= \sum_{i,j=1}^ng^{ij} \nabla_i \nabla_j:C^{\infty}(M)\to  C^{\infty}(M)$ 
is the laplacian  of $g$}.   
$
Then, for a certain constant $C$ we have  $\Delta_g (\lambda +C) = -2(n+1) (\lambda+C)$. Thus, $\lambda+C$ is an eigenfunction of $\Delta_g$ with negative  eigenvalue $-2(n+1)$. Since 
the metric $g$ is negative-definite and the manifold is closed, 
  laplacian  of $g$ is positive definite on nonconstant functions  implying  $\lambda+C\equiv \const$. Thus, 
  $\lambda$ is constant.   Theorem \ref{main2} is proved.  \\[.1cm]

{\bf Proof of Corollaries \ref{Cor1}, \ref{Cor2}.}
It is well-known that if a  manifold $(\hat M, \hat g)$  is decomposable, then  there exists a symmetric tensor  $\hat a= (\hat a_{ij})$, $i,j=0,...,n$ such that $\hat a\ne \const\cdot  \hat g$ for every  $\const \in \mathbb{R}$ and such that    its covariant derivative vanishes: $\hat \nabla_k \hat a_{ij}\equiv 0$. We denote by $\mu$ the $(0,0)-$componenent of $\hat a$,  by $\lambda_i$ the $(0,i)-$component of    $\hat a$ (the symmetric $(i,0)-$component  is also $\lambda_i$), and by $a_{ij}$ the $(i,j)-$component of $\hat a$ for  $i,j=1,...,n$, so  that  the matrix  of $\hat a$ is {
\begin{equation}\label{hata} 
(\hat a_{ij})= \begin{pmatrix} 
\mu & \lambda_1& \dots & \lambda_n \\
\lambda_1&  a_{11}   &  \dots          &    a_{1n}       \\ 
\vdots &     \vdots   &   &         \vdots    \\ 
\lambda_n&    a_{n1}  &  \dots       & a_{nn} 
\end{pmatrix}
\end{equation}
}

The components of $\mu, \lambda_i, a_{ij}$ can a priori depend on $t$. 
For  a fixed $t$ (say, for $t=1$), one can view $\mu, \lambda_i, a_{ij}$ as geometrical objects on $M$: 
$\mu$ is a function on $M$, $\lambda_i$ is an $(0,1)-$tensor on $M$, and $a_{ij}$ is a symmetric $(0,2)-$tensor on $M$ (i.e., if we change the local coordinate system on $M$ the componenents of $\lambda_i$ and $a_{ij}$ change according to the tensor rules). We will denote by $\nabla$ ($\hat \nabla$, resp.)  the covariant derivative in the sense of $g$ ($\hat g$, resp.) 
and by $\Gamma_{ij}^k$  ($\hat \Gamma_{ij}^k$, resp.) the corresponding  Christoffel symbols. 
 We will need the following   
\begin{Lemma} \label{Lem1} Let $\hat a$ given by \eqref{hata} satisfy $\hat \nabla \hat a=0$. 
Then, the tensors $\lambda_i, $ $a_{ij}$, and the function $\mu$  on $M$ satisfy  (we assume $t=1$)
\begin{eqnarray} \label{basic} 
\nabla_ka_{ij} &=&- \lambda_i g_{jk} -\lambda_j g_{ik}, \\ 
\nabla_j \lambda_i & = & a_{ij}- \mu g_{ij}, \label{vnb} \\
\nabla_i\mu &=& 2 \lambda_i. \label{mu}  
\end{eqnarray} 
\end{Lemma}  
{\bf Proof.} Let us calculate  $\hat \Gamma^{i}_{jk}$   in terms of  $g_{ij}$ and  $\Gamma^{i}_{jk}$ at the point $(1, x_1,...,x_n)$ of  the cone $\hat M$: substututing \eqref{hatg} in  $\hat \Gamma^{i}_{jk}= \tfrac{1}{2} \sum_{h = 0}^n \hat g^{i h }\left(\partial_k \hat g_{j h }+  \partial_j \hat g_{h  k} - \partial_h \hat g_{jk} \right) $  we obtain
  \begin{equation} \label{gamma} 
 \begin{array}{lrclr} \hat \Gamma_{j0}^0 = \hat \Gamma_{0j}^0= 0 & \forall j\in 0,...,n  &|&  \hat\Gamma^0_{jk}=-  g_{jk} & \forall j,k \in 1,...,n\\
\hat \Gamma_{j0}^j = \hat \Gamma_{0j}^j= 1 & \forall j\in 1,...,n &|& 
\hat \Gamma_{j0}^i = \hat \Gamma_{0j}^i= 0 & \forall i\ne j\in 1,...,n.  \\
\hat\Gamma^i_{jk}= \Gamma^i_{jk} & \forall i,j,k \in 1,...,n &|&&
\end{array}  
\end{equation} 
Substituting  $\eqref{hata}$ and  $\eqref{gamma}$   in the equation $\hat\nabla_k\hat a_{ij}=0$, we obtain that  for every $i,j,k\in 1,...,n$  
$$
0=\hat\nabla_k\hat a_{ij}= \partial_k a_{ij} -\hat\Gamma_{kj}^0 \hat a_{i0} -  \hat\Gamma_{ik}^0 \hat a_{0j}- \sum_{h=1}^n 
\left[\hat\Gamma_{kj}^h \hat a_{ih} +  \hat\Gamma_{ik}^h \hat a_{hj}\right]  =\nabla_k a_{ij} + g_{kj}\lambda_i + g_{ik} \lambda_j, 
$$
which proves   \eqref{basic}. Similarly, substituting 
 $\eqref{hata}$ and  $\eqref{gamma}$   in   $\hat\nabla_j\hat a_{i0}=0$  we obtain \eqref{vnb}, and substituting $\eqref{hata}$ and  $\eqref{gamma}$ in $\hat\nabla_i\hat a_{00}=0$ we obtain \eqref{mu}. Lemma \ref{Lem1} is proved.

\begin{Lemma} \label{cor2} 
The $(0,1)-$tensor $\lambda_i$ is the  differential of a certain function $\lambda$ on $M$,
 i.e., $\lambda_i= \nabla_i \lambda = \partial_i\lambda$.  Moreover, the function 
 $\lambda$ satisfies the equation \eqref{tanno}.  
  Moreover, if $\lambda$ is constant, then $\hat a$  is proportional to $\hat g$ (with a constant coefficient of proportionality).  
\end{Lemma} 
{\bf Proof.} We multiply (\ref{basic}) by $g^{ij}$ (which is the 
 dual tensor to $g_{ij}$: $\sum_{h=1}^n g^{ih}g_{hj}= \delta_j^i$) 
 and sum over $i$ and $j$: since $\nabla_k g^{ij}=0$  we obtain $\nabla_{k} \sum_{i,j=1}^n a_{ij} g^{ij} = -2 \lambda_k$. Thus, 
 $\lambda_k=  \nabla_{k}\left( -\tfrac{1}{2} \sum_{i,j=1}^n a_{ij} g^{ij}\right)= \nabla_k \lambda$ for the function $\lambda:= -\tfrac{1}{2} \sum_{i,j=1}^n a_{ij} g^{ij}$. 
 Now, covariantly differentiating \eqref{vnb}, replacing $\lambda_i$ by $\nabla_i\lambda $ 
 and replacing  the covariant derivatives of $a_{ij}$ and $\mu$ using \eqref{basic} and \eqref{mu} we obtain 
 $$\begin{array}{rr}    
0= &  \nabla_k\left(\nabla_j\lambda_i  -a_{ij}+ \mu g_{ij}\right)= \nabla_k\nabla_j\nabla_i \lambda - \nabla_k a_{ij} + \nabla_k \mu \cdot g_{ij} \\ = &  \nabla_k\nabla_j\nabla_i \lambda  +
   \nabla_i \lambda\cdot  g_{jk} + \nabla_j \lambda \cdot  g_{ik} + 2 \nabla_k \lambda \cdot g_{ij},  \end{array}
$$ which is the equation \eqref{tanno}. 

If $\lambda$ is constant,   $\mu$  is constant by \eqref{mu}. Then,   
 \eqref{vnb} implies  $a=\mu\cdot g$. Since $\lambda_i=\partial_i\lambda=0$, we have     $\hat a = \mu \cdot \hat g$, i.e., $\hat a $ is proportional to $\hat g$ at every point of $\hat M$ with $t=1$. Since $\hat a$ and $\hat g$ are covariantly constant, $\hat a $ is proportional to $\hat g$ at every point of $\hat M$.   Lemma \ref{cor2} is proved. \\[.1cm] 
 
\begin{Rem} 
 Corollaries \ref{Cor1}, \ref{Cor2} easily follow from Theorems \ref{main1}, \ref{main2} and Lemma  \ref{cor2}. \end{Rem}

{\bf Certain generalizations.}  One can easily generalize our proof of Theorem \ref{main1} for higher Gallot equations $E_p$  introduced in \cite[Section 4]{G}: for every $p\in \mathbb{N}$ the equation $E_p$ is 
\begin{equation} \label{Ep}
D^{p+1} f(Y_1,...,Y_{p+1}) + \sum_{1\le s\le \tfrac{p+1}{2} } \sum_{\sigma\in S_{p+1}} \lambda(s,\sigma) \left(  g^s\otimes D^{p+1-2s} f \right)(Y_{\sigma(1)},...,Y_{\sigma(p+1)})=0,
\end{equation}   where $f$ is the unknown function,  $S_{p+1} $ denotes the set of all permutations of $\{1,...,p+1\}$,  $\lambda(s,\sigma)$ denotes certain numbers  depending on $s\in 1,..., [\tfrac{p+1}{2}]$ and on $\sigma\in S_{p+1}$ whose precise values are  not important for our proof,  $Y_1,...,Y_{p+1}$ are arbitrary vector fields, and $D^{k}$ denotes the $k-$th covariant derivative (so for example $D^2f (X,Y)= \sum_{i,j=1}^n X^i Y^j \nabla_j \nabla_i f$).

\begin{Th} \label{main3}
Let $g$ be  
a light-line-complete 
connected  pseudo-Riemannian metric of indefinite signature on a closed   $n-$dimensional manifold $M^n$. Then, every solution of \eqref{Ep} is constant.  
\end{Th} 
{\bf Proof. } We take a light-line geodesic $\gamma$, 
 and take arbitrary vector fields $Y_i$ such that at every point of the geodesic $\gamma$ we have 
  $Y_i= \dot\gamma$.  Since $g(\dot\gamma,\dot\gamma)= 0$  and $s\ge 1$, we obtain  \\
  $\left(  g^s\otimes D^{p+1-2s} f \right)(Y_{\sigma(1)},...,Y_{\sigma(p+1)})=0$. Then,   $D^{p+1}( \dot\gamma,...,\dot\gamma)=0$ implying $\tfrac{d^{p+1}}{dt^{p+1}}f(\gamma(t))\equiv 0$ implying 
 $f= \const_p t^p + ... + \const_0$. Since the manifold is compact, the function $f$  must be bounded implying $\const_p=...=\const_1=0$. Thus, the function $f$ must be constant along every light-line geodesic.  Since every two points of a connected pseudo-Riemannian manifold of indefinite signature  can be connected by a sequence of light-line geodesics, the function $\lambda$ is   a constant. Theorem  \ref{main3} is proved.

 Another possible generalization is due to the observation that in our proof of Corollaries 
 \ref{Cor1}, \ref{Cor2} we actually used the existence of a covariantly-constant symmetric 
 $(0,2)-$tensor $\hat a_{ij}\ne \const \cdot \hat g_{ij}$ only. Decomposability of  the metric $\hat g$ implies   the existence of such a  tensor $\hat  a$,  but not vice versa: 
  in the pseudo-Riemannian case there  exist metrics $g$ 
  admitting covariantly-constant symmetric $a\ne \const \cdot  g$, see \cite{ks}.   
  So, in fact we proved 
  
  \begin{Cor}  \label{Cor3} 
Let $g$ be  
a light-line-complete 
 pseudo-Riemannian metric of indefinite signature   on a closed  $n-$dimensional manifold $M^n$. 
 Then,   every  symmetric $(0,2)-$tensor $\hat a_{ij}$ on the corresponding cone $(\hat M, \hat g) $  such that 
 $\hat \nabla_k \hat a_{ij}\equiv 0$ is proportional to $\hat g_{ij}$.
 
\end{Cor}

\begin{Cor} \label{Cor4} 
Let $g$ be  
a complete  negative-definite 
 pseudo-Riemannian metric   on a closed   $n-$dimensional manifold $M^n$. 
Then,   every  symmetric $(0,2)-$tensor $\hat a_{ij}$ on the  corresponding  cone $(\hat M, \hat g) $   such that 
 $\hat \nabla_k \hat a_{ij}\equiv 0$ is proportional to $\hat g_{ij}$.
\end{Cor}

{{\bf Acknowledgement:} } 
We  thank    Deutsche Forschungsgemeinschaft (Priority Program 1154 --- Global Differential Geometry)   and FSU Jena for partial financial support, and Vicente Cortes,  Thomas Leistner and Dmitri Alekseevsky   for useful discussions.


\begin{thebibliography}{99}

\bibitem{ACGG} D. V. Alekseevsky, V. Cortes, A. S. Galaev and T. Leistner, {\it   Cones over pseudo-Riemannian manifolds and their holonomy,} to appear in Journal f\"ur die reine und angewandte Mathematik (Crelle's Journal), arXiv:0707.3063. 

\bibitem{couty}  R. Couty, \emph{Transformations
 infinit$\acute{e}$simales projectives,}
  C. R. Acad. Sci. Paris {\bf 247}(1958), 804--806, MR0110994, Zbl 0082.15302.



\bibitem{G} S. Gallot, {\it \'Equations diff\'erentielles caract\'eristiques de la sph\`ere,}  Ann. Sci. \'Ecole
Norm. Sup. (4), {\bf 12}(2), 235--267, 1979.

\bibitem{hall} G. S. Hall,  D. P.        Lonie,  {\it Projective Structure of Einstein Spaces in General
Relativity,} accepted to Classical and Quantum Gravity.

\bibitem{einstein} 
V. Kiosak,  V. S. Matveev,  \emph{ Complete Einstein metrics are geodesically rigid,} 
 Comm. Math. Phys. {\bf 289}(1),   383--400,  2009, 
  	arXiv:0806.3169.

\bibitem{KM} V. Kiosak, 
V. S. Matveev, 
{\em  Proof of projective Lichnerowicz conjecture for pseudo-Riemannian metrics with degree of mobility greater than two,}   arXiv:0810.0994 

\bibitem{ks} G. I. Kruzhkovich, A. S. Solodovnikov, {\it 
Constant symmetric tensors in Riemannian spaces,} 
Izv. Vyssh. Uchebn. Zaved. Matematika {\bf 1959}  no. 3 (10), 147--158. 

\bibitem{Obata} M.  Obata, \emph{
Riemannian manifolds admitting a solution of a certain system of differential equations,}  Proc. U.S.-Japan Seminar in Differential Geometry (Kyoto, 1965) pp. 101--114.  

\bibitem{sol} A. S. Solodovnikov, {\it Projective transformations of Riemannian spaces,}
 Uspehi Mat. Nauk (N.S.) {\bf 11}(1956), no. 4(70), 45--116, MR0084826, Zbl 0071.15202.
 
\bibitem{Tanno} S.  Tanno, \emph{Some differential equations on Riemannian manifolds,} 
J. Math. Soc. Japan {\bf 30}(1978), no. 3, 509--531.



\bibitem{Vries} H. L. de  Vries, {\it \"Uber Riemannsche R\"aume, die
infinitesimale konforme Transformationen gestatten,}
 Math. Z. {\bf 60}(1954), 328--347, MR0063725, Zbl 0056.15203.


\end{thebibliography}
\end{document}